\documentclass[12pt]{article} 
\usepackage{epsfig}
\usepackage{psfrag}
\usepackage{amsbsy}
\usepackage{latexsym}
\usepackage[mathscr]{eucal}
\usepackage{amsfonts}
\usepackage{amssymb}

\def\inpro#1{\langle #1\rangle}
\def\tilpsi{{\widetilde \psi}}

\def\eps{\epsilon}
\def\calD{{\cal D}}
\def\calV{{\cal V}}
\def\v{v}
\def\iv{I_\v}
\def\ivp{I_{\v'}}
\def\bivp{{\bar I}_{\v'}}
\def\bivph{{\bar I}_{\v{\vphantom '}}}
\def\biv{{\bar I}_{\v}}
\def\lam{\lambda}
\def\sig{\sigma}
\def\Ome{\Omega}
\def\del{\delta}
\def\Lam{\Lambda}
\def\tilh{{\widetilde h}}
\def\k{\kappa}

\begin{document}

 \title{Approximation using scattered shifts of a multivariate 
function \thanks{ 
This work has been supported  by the Office of Naval 
Research Contracts   ONR-N00014-03-1-0051, ONR/DEPSCoR N00014-03-1-0675,  
ONR/DEPSCoR N00014-05-1-0715;
the Army Research Office 
Contracts  DAAD 19-02-1-0028,   W911NF-05-1-0227, and   W911NF-07-1-0185;  the National Institute of General Medical 
Sciences under Grant NIH-1-R01-GM072000-01;  the National Science 
Foundation  under Grants DMS-0221642,  DMS-9872890, DMS-354707, DBI-9983114,   
ANI-0085984 and DMS-0602837}}
\author{Ronald DeVore and Amos Ron}
\hbadness=10000
\vbadness=10000
\newtheorem{lemma}{Lemma}[section]
\newtheorem{prop}[lemma]{Proposition}
\newtheorem{cor}[lemma]{Corollary}
\newtheorem{theorem}[lemma]{Theorem}
\newtheorem{remark}[lemma]{Remark}
\newtheorem{example}[lemma]{Example}
\newtheorem{definition}[lemma]{Definition}
\newtheorem{proper}[lemma]{Properties}
%
\def\vp{\varphi}
\def\<{\langle}
\def\>{\rangle}
\def\t{\tilde}
\def\i{\infty}
\def\e{\varepsilon}
\def\sm{\setminus}
\def\nl{\newline}
\def\Chi{\raise .3ex \hbox{\large $\chi$}} 
\def\lsima{\hbox{\kern -.6em\raisebox{-1ex}{$~\stackrel{\textstyle<}{\sim}~$}}\kern -.4em}
\def\lsim{\hbox{\kern -.2em\raisebox{-1ex}{$~\stackrel{\textstyle<}{\sim}~$}}\kern -.2em}\def\[{\Bigl [}
\def\]{\Bigr ]}
\def\({\Bigl (}
\def\){\Bigr )}
\def\[{\Bigl [}
\def\]{\Bigr ]}
\def\({\Bigl (}
\def\){\Bigr )}
\def\L{£}
\def\pr{{\rm Prob}}
\newcommand{\iref}[1]{(\ref{#1})}
\def\ds{\displaystyle}
\def\ev#1{\vec{#1}}     
\newcommand{\lt}{\ell_{2}(\nabla)}
\def\Supp#1{{\rm supp\,}{#1}}
\def\e{\epsilon}
\def\R{\mathbb{R}}
\def\E{\mathbb{E}}
\def\nl{\newline}
\def\T{{\relax\ifmmode I\!\!\hspace{-1pt}T\else$I\!\!\hspace{-1pt}T$\fi}}
\def\N{\mathbb{N}}
\def\Z{\mathbb{Z}}
\def\N{\mathbb{N}}
\def\Zd{\Z^d}
\def\Q{\mathbb{Q}}
\def\C{\mathbb{C}}
\def\Rd{\R^d}
\def\gsim{\mathrel{\raisebox{-4pt}{$\stackrel{\textstyle>}{\sim}$}}}
\def\sime{\raisebox{0ex}{$~\stackrel{\textstyle\sim}{=}~$}}
\def\lsim{\raisebox{-1ex}{$~\stackrel{\textstyle<}{\sim}~$}}
\def\div{\mbox{ div }}
\def\M{M}  \def\NN{N}                  
\def\L{{\ell}}               
\def\Le{{\ell_1}}            
\def\Lz{{\ell_2}}
\def\Let{{\tilde\ell_1}}     
\def\Lzt{{\tilde\ell_2}}
\def\Ltw{\ell_\tau^w(\nabla)}
\def\t#1{\tilde{#1}}
\def\la{\lambda}
\def\La{\Lambda}
\def\ga{\gamma}
\def\BV{{\rm BV}}
\def\Ga{\eta}
\def\al{\alpha}
\def\cZ{{\cal Z}}
\def\argmin{\mathop{\rm argmin}}
\def\argmax{\mathop{\rm argmax}}

\def\prob{\mathop{\rm prob}}

\def\cO{{\cal O}}
\def\cA{{\cal A}}
\def\cC{{\cal C}}
\def\cF{{\cal F}}
\def\bu{{\bf u}}
\def\bz{{\bf z}}
\def\bZ{{\bf Z}}
\def\bI{{\bf I}}
\def\cE{{\cal E}}
\def\cD{{\cal D}}
\def\cG{{\cal G}}
\def\cI{{\cal I}}
\def\cJ{{\cal J}}
\def\cM{{\cal M}}
\def\cN{{\cal N}}
\def\cT{{\cal T}}
\def\cU{{\cal U}}
\def\cV{{\cal V}}
\def\cW{{\cal W}}
\def\cL{{\cal L}}
\def\cB{{\cal B}}
\def\cG{{\cal G}}
\def\cK{{\cal K}}
\def\cS{{\cal S}}
\def\cP{{\cal P}}
\def\cQ{{\cal Q}}
\def\cR{{\cal R}}
\def\cU{{\cal U}}
\def\bL{{\bf L}}
\def\bK{{\bf K}}
\def\bC{{\bf C}}
\def\X{X\in\{L,R\}}
\def\ph{{\varphi}}
\def\D{{\Delta}}
\def\H{{\cal H}}
\def\bM{{\bf M}}
\def\bx{{\bf x}}
\def\bG{{\bf G}}
\def\bP{{\bf P}}
\def\bW{{\bf W}}
\def\bT{{\bf T}}
\def\bV{{\bf V}}
\def\bv{{\bf v}}
\def\bt{{\bf t}}
\def\bz{{\bf z}}
\def\bw{{\bf w}}
\def \span{{\rm span}}
\def \meas {{\rm meas}}
\def\rhom{{\rho^m}}
\def\lll{\langle}
\def\argmin{\mathop{\rm argmin}}
\def\argmax{\mathop{\rm argmax}}
\def\dJ{\nabla}
\newcommand{\ba}{{\bf a}}
\newcommand{\bb}{{\bf b}}
\newcommand{\bc}{{\bf c}}
\newcommand{\bd}{{\bf d}}
\newcommand{\bs}{{\bf s}}
\newcommand{\bff}{{\bf f}}
\newcommand{\bp}{{\bf p}}
\newcommand{\bg}{{\bf g}}
\newcommand{\by}{{\bf y}}
\newcommand{\br}{{\bf r}}
\newcommand{\be}{\begin{equation}}
\newcommand{\ee}{\end{equation}}
\newcommand{\bea}{$$ \begin{array}{lll}}
\newcommand{\eea}{\end{array} $$}
\def \Vol{\mathop{\rm  Vol}}
\def \mes{\mathop{\rm mes}}
\def \Prob{\mathop{\rm  Prob}}
\def \exp{\mathop{\rm    exp}}
\def \sign{\mathop{\rm   sign}}
\def \vphi{{\varphi}}
\def \csp{\overline \mathop{\rm   span}}

%
\newcommand{\beqn}{\begin{equation}}
\newcommand{\eeqn}{\end{equation}}
\def\beginproof{\noindent{\bf Proof:}~ }
\def\endproof{\hfill\rule{1.5mm}{1.5mm}\\[2mm]}

\newenvironment{Proof}{\noindent{\bf Proof:}\quad}{\endproof}

\renewcommand{\theequation}{\thesection.\arabic{equation}}
\renewcommand{\thefigure}{\thesection.\arabic{figure}}

\makeatletter
\@addtoreset{equation}{section}
\makeatother

\newcommand\abs[1]{\left|#1\right|}
\newcommand\clos{\mathop{\rm clos}\nolimits}
\newcommand\trunc{\mathop{\rm trunc}\nolimits}
\renewcommand\d{d}
\newcommand\dd{d}
\newcommand\diag{\mathop{\rm diag}}
\newcommand\dist{\mathop{\rm dist}}
\newcommand\diam{\mathop{\rm diam}}
\newcommand\cond{\mathop{\rm cond}\nolimits}
\newcommand\eref[1]{(\ref{#1})}
\newcommand\Hnorm[1]{\norm{#1}_{H^s([0,1])}}
\def\int{\intop\limits}
\renewcommand\labelenumi{(\roman{enumi})}
\newcommand\lnorm[1]{\norm{#1}_{\ell_2(\Z)}}
\newcommand\Lnorm[1]{\norm{#1}_{L_2([0,1])}}
\newcommand\LR{{L_2(\R)}}
\newcommand\LRnorm[1]{\norm{#1}_\LR}
\newcommand\Matrix[2]{\hphantom{#1}_#2#1}
\newcommand\norm[1]{\left\|#1\right\|}
\newcommand\ogauss[1]{\left\lceil#1\right\rceil}
\newcommand{\QED}{\hfill
\raisebox{-2pt}{\rule{5.6pt}{8pt}\rule{4pt}{0pt}}%
  \smallskip\par}
\newcommand\Rscalar[1]{\scalar{#1}_\R}
\newcommand\scalar[1]{\left(#1\right)}
\newcommand\Scalar[1]{\scalar{#1}_{[0,1]}}
\newcommand\Span{\mathop{\rm span}}
\newcommand\supp{\mathop{\rm supp}}
\newcommand\ugauss[1]{\left\lfloor#1\right\rfloor}
\newcommand\with{\, : \,}
\newcommand\Null{{\bf 0}}
\newcommand\bA{{\bf A}}
\newcommand\bB{{\bf B}}
\newcommand\bR{{\bf R}}
\newcommand\bD{{\bf D}}
\newcommand\bE{{\bf E}}
\newcommand\bF{{\bf F}}
\newcommand\bH{{\bf H}}
\newcommand\bU{{\bf U}}
\newcommand\cH{{\cal H}}
\newcommand\sinc{{\rm sinc}}
\def\enorm#1{| \! | \! | #1 | \! | \! |}

\newcommand{\dm}{\frac{d-1}{d}}

\let\bm\bf
\newcommand{\bbeta}{{\mbox{\boldmath$\beta$}}}
\newcommand{\bal}{{\mbox{\boldmath$\alpha$}}}
\newcommand{\bbi}{{\bm i}}

\newif\ifNZB
\newcommand\NZB[1]{\ifNZB \marginpar{\raggedright \scriptsize NZB:\\#1}
 \else \fi}
\newcommand{\FText}[1]{\mbox{#1}}
\makeatletter
\newcommand{\tr}{{\mathop{\operator@font T}\nolimits}}
\newcommand{\mod}{\mathop{\operator@font mod}}
\makeatother
\newcommand{\UArrow}[4]{
 \begin{array}{ll}
  #1&\stackrel{#2}\longrightarrow\\
  #3&\;\raisebox{1ex}{$\nearrow$}\mkern-14mu_{#4}
 \end{array}}
\newcommand{\DArrow}[4]{
 \begin{array}{ll}
  \stackrel{#1}\longrightarrow&#2\\
  \mkern-10mu_{#3}\mkern-22mu\raisebox{1ex}{$\searrow$}&#4
 \end{array}}
\newcommand{\fig}[3]{\par\begin{figure}[ht]
  \centerline{\epsfbox{#1.eps}}\caption{#3}\label{fig#2}\end{figure}}
\newcommand{\dI}{\Delta}
\maketitle
\date{}
\begin{abstract}
The approximation of a general $d$-variate
function $f$ by the shifts $\phi(\cdot-\xi)$,
$\xi\in\Xi\subset \Rd$, of a fixed function $\phi$  occurs in many
applications such as data fitting, neural networks, and learning theory. 
When $\Xi=h\Z^d$
is a dilate of the integer lattice, 
there is a rather
complete understanding of the approximation problem \cite{BDR,Johnson1} using
Fourier techniques.  However, in most applications the {\it center} set $\Xi$
is either given, or can be chosen with complete freedom.  In both of these cases,  the shift-invariant setting is too restrictive.  This paper 
studies the approximation problem in the case $\Xi$ is arbitrary.  It 
establishes approximation theorems whose error bounds reflect the local 
density of the points in $\Xi$. Two different settings are analyzed.
The first is when the set $\Xi$ is prescribed in advance. In this case, the 
theorems of this paper show that, in analogy with the classical
univariate spline approximation,  improved approximation occurs in 
regions where the density is high.  The second setting corresponds to the
problem of non-linear approximation.  In that setting the set $\Xi$ can 
be chosen using information about the 
target function $f$.  
We discuss how to `best' make these choices and  give estimates for the 
approximation error.  \\

\medskip
\noindent
{\bf AMS subject classification:} 42C40, 46B70, 26B35, 42B25\\

\noindent
{\bf Key Words:}
image/signal processing, computation, nonlinear approximation, optimal
approximation, radial basis functions, scattered data, thin-plate splines,
surface splines, approximation order
\end{abstract}
 %
\section{Introduction}
\label{intro}

The mathematical problem of data fitting in the $d$-variate Euclidean space
$\Rd$ has vast applications in science and engineering. Many
algorithms  
address this problem by approximating the data by a linear combination
$F=\sum_{\xi\in\Xi}c(\xi)\phi(\cdot-\xi)$, with $\Xi\subset\Rd$, and
$\phi$ a carefully-selected, often radial, function
defined on $\Rd$.  One of the primary motivations for this    
approach  is that
if the data themselves are defined on $\Xi$ and $\phi$ is chosen
properly, then there is a  unique
function related to the above that  interpolates the data
\cite{Micchelli,scho1,scho2}.
For example, if $\phi$ is the so-called surface spline
(a fundamental solution of the  $m$-fold iterated Laplacian) and $\Xi$ is a 
given finite set of points in $\Rd$, then for given data $(\xi,y_\xi)$, 
$\xi\in\Xi$, there is a unique   interpolant to the data from the span
$S_\Xi(\phi)$ of the $\phi(\cdot-\xi)$, $\xi\in\Xi$.%
\footnote{The interpolant is actually selected
from  a space of the form $S'_\Xi(\phi)\oplus P$, with $S'_\Xi(\phi)$
a certain subspace of $S_\Xi(\phi)$, and $P$ a finite dimensional space
of polynomials, that depends only on $\phi$.}
We refer to the book \cite{Wendland}
for more details on radial basis functions in general, and their
use in interpolation, in particular.

The problem of estimating the interpolation error in the above setting was
studied extensively in the literature.  We refer the reader to 
\cite{MaNel,WuSchaback}, where the interpolated function
$f$ is assumed to come from the so-called ``native space''
(this approach originated in the work of Duchon, \cite{Duchon1,Duchon2})
and to the more recent \cite{Schaback,YoonJAT01,YoonSIMA01,NW}, 
where error estimates are established for general smooth functions in 
Sobolev spaces.  It should be noted 
that the interpolation problem is usually analyzed for functions defined on 
bounded subdomains of $\Rd$, and it is well-understood that the 
interpolation error for this setup suffers significantly from the so-called 
``boundary effect''. Typically, the rates of decay of the error for 
smooth functions are about half the corresponding decay rates that are valid 
in the boundary-free shift-invariant case, \cite{Johnson2,Hangel}.

Interpolation is not necessarily the best approach to the data fitting 
problem for various reasons including possible noise in the data, the
computational overhead, and possible lack of stability in the algorithms.
If  the data are given by a function $f$ defined on $\Rd$, i.e.  
$y_\xi=f(\xi)$  (or $y_\xi\approx f(\xi)$ in the noisy version of the 
problem), the primary question is how well can $f$ be approximated in a 
given metric (typically $L_p$-norms) from the given information.  This is
governed, and in part determined, by the  related question of how well $f$ 
can be approximated from the span of the translates
$\phi(\cdot-\xi)$, $\xi\in \Xi$, in the given metric.

In this paper, we shall be solely concerned with  the latter {\it 
approximation} problem.   We start with a countable set $\Xi$  of points in 
$\Rd$ and define $S_\Xi(\phi)$  to be the set of all functions which are 
finite linear combinations of the shifts $\phi(\cdot-\xi)$, $\xi\in \Xi$.  We 
are interested in how well a given function $f\in L_p(\Rd)$ can be 
approximated (in the $L_p$-norm)  by the elements of $S_\Xi(\phi)$ 
(more precisely by elements in  the closure of this space in the given 
metric).  Such approximation problems have been well studied, especially in 
the case that $\Xi$ is a dilate of the integer lattice ($\Xi=h\Zd$ with 
$h>0$), \cite{BRquasi92,Buhmann,BR92,BDR,Ron92NMAT,Johnson1}.  
The case where the centers $\Xi$ are scattered was studied in 
\cite{BDL,DyRo,BuhRonproc92,Johnson3,YoonCA01}. The error bounds in all
these references are given in terms of a global mesh density parameter.
In contrast, error bounds that depend on the local density of the
scattered centers (i.e., provide improved error bounds on subdomains
that contain dense clusters of centers) are less studied and less understood,   
even though it is often the natural setting in applications.  
The most notable exception, is, of course, 
spline approximation in one variable. The fact that the error bounds in linear
approximation by splines reflects the local mesh ratio, \cite{deboor_book},
is a key property of spline approximation. Furthermore, the development
and analysis of non-linear approximation schemes for univariate splines
\cite{Pencho,DePop,DeYu} presented the first challenge for the
development of the substantial theory of non-linear approximation.
In more than one 
variable, however, far less is known.
We refer to \cite{Powell}, where low-rate strongly local error
estimates are established, and to the approximation scheme based on  the
``power function approach'' in \cite{WuSchaback}.  

We shall consider two types of problems for scattered center 
approximation.  In the first, we assume that the set $\Xi$ is fixed and we 
derive results that show improved approximation in regions where the density 
is high.  These results are described in \S \ref{sect3} and \S 
\ref{seclower}.  The second setting that we consider allows the centers to 
be chosen dependent on the function $f$. The basic goal is to establish 
error bounds that depend on the cardinality $N$ of the chosen center set 
$\Xi$.  This is a form of nonlinear approximation known as $N$-term 
approximation which has been well studied in other settings, primarily for 
wavelet bases.  Our result  here is similar to the results on nonlinear 
wavelet approximation.  We show that a function can be approximated in 
$L_p(\Rd)$ with error $O(N^{-s/d})$ once it lies in the Triebel-Lizorkin space 
$F^s_{\tau,q}(\Rd))$  where $s$, $p$, and $\tau$ are related (as in the 
Sobolev embedding theorem)  by $\frac{1}{\tau}-\frac{1}{p}=\frac{s}{d}$ and 
$q=(1+\frac{s}{d})^{-1}$. From this result and standard embeddings for  
Triebel-Lizorkin spaces, we derive corresponding theorems for $N$-term 
approximation in terms of the Besov classes.   While our actual results
in this direction are close in nature to the wavelet results, the 
non-linear approximation algorithm that leads to the above error bounds 
differs from its wavelet counterpart: the thresholding algorithm that is 
employed in the wavelet case is sub-optimal in the present case; as such, we 
introduce and analyse a more sophisticated algorithm.  Details of this 
result are given in \S \ref{secnonlinear}. 

We begin in the following section by describing the assumptions we make about 
$\phi$ and $\Xi$.  We then give our first results for the linear approximation
problem in \S \ref{sect3}.
In \S \ref{secwave} we recall some results about 
wavelet decompositions and the use of such decompositions in the 
characterization of smoothness spaces (Triebel-Lizorkin spaces;
they include the more standard Sobolev spaces).  We stress
that our paper is not concerned with wavelets: we merely use
wavelets as a tool for defining our approximation schemes.
In \S \ref{seclower}, we
complete our  study of the linear approximation problem.  Finally,
we prove in \S \ref{secnonlinear} our results on nonlinear approximation.

We shall treat approximation on all of $\Rd$.  In most applications, one 
would be interested in the case of approximation on a compact domain 
$\Omega$ of $\Rd$.  Results on domains can be derived easily from our 
results (if the approximand $f$ is defined on a domain that includes
$\Omega$ in its interior, and if one agrees to allow centers outside
$\Omega$) but we do not pursue this here.  


\section{The setting}
\label{sect2}
We describe in this section the setting that will
be analyzed in the first part of the paper.  There are two main ingredients 
in our setting.  The first is the set $\Xi$ of centers which can be 
allowed.  We do not make direct assumptions on the geometry of the
set $\Xi$: almost any set $\Xi$ will do. Once the set $\Xi$ is given,
we associate it with a density function
$$h: \Rd\to \R_+.$$
The value $h(t)$ of the
density function depends strongly on the local density of $\Xi$ around
$t$: roughly speaking, there should be $L$ centers of $\Xi$ is  a  ball
of radius $L'h(t)$ centered at $t$, with $L,L'$ dependent only on
$\phi$ and some parameters that we choose for our
approximation scheme (and that we fix throughout).
The density $h(t)$ depends also on the geometry
of the centers around $t$.
This dependence is generally mild; also, our assumptions
never spell out this dependence explicitly. It is embedded implicitly in
other assumptions.

\medskip\noindent
We assume that our set $\Xi$ of centers 
is a countable set in  $\Rd$ and 
\vskip .15in
{\bf A1:} {\it {\rm (i)} Any finite
ball in $\Rd$ contains a finite number of points from $\Xi$.   

\ \ \ \ \ \  {\rm (ii)}   For each integer $n$, there is an $R=R(n)$ such that each ball of radius $R$

\ \ \ \ \ \ \ \ \ \ contains at least $n$ points from $\Xi$.}
\vskip .15in

\noindent
Property (i) prevents the occurrence of
accumulation points in the set $\Xi$ while property
(ii) prohibits the existence of arbitrarily large regions on which there are no 
centers.  Neither of the two conditions in {\bf A1} is essential, and the entire
assumption is adopted merely in order to simplify the presentation and the 
analysis.
  
The first {\it essential} ingredient in our setting is the identification of
the class of basis functions $\phi$ that our analysis applies to.
We will make assumptions on $\phi$ that are 
convenient and make the ensuing analysis as transparent as possible and yet 
general enough to be valid for certain (but not all) types
of multivariate splines and radial basis functions.

We assume that $\phi$ is a
locally integrable function which, when viewed on all of $\Rd$, is a
tempered distribution.  
We denote by 
$$S_\Xi(\phi)$$
the   finite linear
span of the functions $\phi(\cdot-\xi)$, $\xi\in\Xi$.
We put forward three assumptions: one on $\phi$, one on $\Xi$, and
one that connects between $\phi$ and $\Xi$.
At the end of this section, we analyse these assumptions for specific 
choices of $\phi$.

Let $C_0^\infty:=C^\infty_0(\Rd)$ denote the set of all functions that  are
infinitely differentiable  with compact support (test functions) and 
$C_0^k$ the analogous spaces of $k$-times continuously
differential functions with compact support ($C_0^{\vphantom k}:=C_0^0)$.
Our first assumption about $\phi$ is:
\vskip .15in
{\bf A2:}  {\it There is a positive integer $\k>0$ and  
a linear operator $T$ mapping $C_0^\k$ into $C_0^{\vphantom \k}$, such that 
for all  $f\in C_0^\k$,  
\beqn
\label{kernel}
f(x)=\int_{\Rd}T(f)(t)\phi(x-t)\,dt, \quad x\in\Rd.
\eeqn
}%
Note that the integration is well-defined: since $\phi$
is locally in $L_1$, and $T(f)\in C_0^{\vphantom k}$,
we have that $T(f)\phi(x-\cdot)\in L_1$, for   
every $x\in \Rd$.

The  typical example  of  $T$ is a homogeneous elliptic differential  
operator of order $\k$ with constant coefficients. In 
this case $\phi$ is its fundamental solution on $\Rd$. Note that
$\hat\phi$ is, in this case,  a smooth function on $\Rd\backslash 0$.
It will be sometimes convenient to add this assumption to {\bf A2}:
$$\hbox{$\hat\phi$ is smooth on $\Rd\backslash 0$}.$$
On the other hand,
assumption {\bf A2} may appear too strong for certain applications,
since it excludes some interesting examples corresponding to 
fractional differentiation.  Those who may wish to extend our theory
in such directions may allow $T(f)$ to have global support; it will still 
need to decay suitably at infinity. All in all, there is some flexibility
in the formulation of {\bf A2}.

That said, {\bf A2}, either in its current variant or in some related
one, is fundamental: it determines the maximal decay
rate of the error that our approximation scheme can yield
(this is the parameter $\k$), and determines the  
space $W(L_p(\Rd),\phi)$ of smooth functions that can be approximated
at this rate. Let us discuss now this latter issue.

The function space  $W(L_p(\Rd),\phi)$ 
corresponds to the Sobolev space $W^\k(L_p(\Rd))$
in the case $T$ is an elliptic  differential operator of order $\k$
(with constant coefficients). For more general $T$, the definition is more
abstract: Fixing $1\le p\le \infty$, we define  the   semi-norm
\beqn
\label{Wseminorm}
|f|_{W(L_p(\Rd),\phi)}:=\|Tf\|_{L_p(\Rd)}, \quad f\in C_0^\infty,
\eeqn
and the norm:
%
\beqn
\label{Wnorm}
\|f\|_{W(L_p(\Rd),\phi)}:=\|f\|_{L_p(\Rd)}+|f|_{W(L_p(\Rd),\phi)}.
\eeqn
Using the above, we define  $W(L_p(\Rd),\phi)$   to be the   completion of  
$C_0^\infty$ in this topology.  Since the norm in \eref{Wnorm}
is stronger than the $L_p$-norm,  $W(L_p,\phi)$ is precisely the 
space of all functions $f\in L_p(\Rd)$ for which there is a sequence 
$(f_n)_{n\ge 1}\subset C_0^\infty$ 
such that $f_n$ converges to $f$ and $(T(f_n))$ converges to a function 
$g\in L_p(\Rd)$ both in the sense of $L_p(\Rd)$.  By making suitable
assumptions on $\phi$ we can conclude that $g$ depends on $f$,
but not on the specific sequence $(f_n)$.
The operator $T$, which was 
initially defined only for test functions,   now   extends naturally to  a 
linear operator on $W(L_p,\phi)$ by defining $Tf:=g$.     

We are guided in the above setup by the following example,
\cite{Duchon1,Duchon2}: if
$T$ is the $m$-fold Laplacian,  then $\k=2m$ and the function $\phi$ is then
the fundamental solution of 
$T$:
\beqn\phi=c_m\cases{|\cdot|^{2m-d},&$d$ odd,\cr
               |\cdot|^{2m-d}\log|\cdot|,&$d$ even.\cr}\label{surface}\eeqn
	       (Here, $|\cdot|$ stands for the Euclidean norm in $\Rd$.)
The function $\phi$ is also   called    a {\it surface spline}.  In this 
case, $W(L_p(\Rd),\phi)$ is simply the Sobolev space $W^{2m}(L_p(\Rd))$ 
equipped with its usual semi-norm and norm.

\medskip
Our remaining  assumption about $\phi$ concerns how well its translates
can be approximated from $S_\Xi(\phi)$.  Consider the translate $
\phi(\cdot  -t)$ where 
$t\in\Rd$ is fixed for the moment.  We look for a local approximation of the
form 
\beqn
\label{approxkernel}
K(\cdot ,t)=\sum_{\xi\in \Xi(t)}A(t,\xi)\phi(\cdot -\xi),
\eeqn
for suitable $t$-dependent coefficients $A(t,\xi)$ and a finite set 
$\Xi(t)\subset \Xi$.  A key component in the success of our approach is
the availability of kernels $K(\cdot,\cdot)$ that are local and bounded
on the one hand, and approximate well the convolution kernel
$(x,t)\mapsto \phi(x-t)$. We break this assumption into two: {\bf A3} deals
with basic qualitative properties of the scheme that is used to define $K$, 
viz., the coefficient functions $A(\cdot,\xi)$. The companion property, {\bf A4},
deals with the way $K$ approximates the convolution kernel.

\vskip .15in
{\bf A3:} {\it  There is an integer   $n'>0$ and a real number $M_0$ such 
that for   any $t\in\Rd$   the  set  $\Xi (t)$ consists of at most $n'$   
points all lying in the ball of radius $M_0$ centered at $t$ and the 
coefficients of the  approximation kernel   {\rm\eref{approxkernel}}  
satisfy   $A(\cdot,\xi)\in L_1(\Rd)$ for all $\xi\in\Xi$. } 
\vskip .15in
\noindent

Similar  to condition {\bf A1}, Condition {\bf A3}   is secondary, and is
formulated and adopted in order to exclude pathological kernels $K$.
This brings us to  
our last assumption. As said, the last assumption
is concerned with the way $K(\cdot,t)$ 
approximates the translate  $\phi(\cdot-t)$. This assumption must be dealt with 
care: the error $E(\cdot,t):=\phi(\cdot-t)-K(\cdot,t)$ should reflect
not only the basic properties of $\phi$, but also the local distribution
of the center set $\Xi$ around $t$. 
To this end, we  define for each $t\in \Rd$
\beqn
\label{defh}
h(t):= \inf \{\rho: A(t,\xi)\equiv 0,\ |\xi-t|\ge \rho\}.
\eeqn
In other words, for each $t\in\Rd$, the only centers from $\Xi$ used in
the approximation kernel $K(\cdot,t)$ lie  in a ball $B_t(h(t))$
of radius $h(t)$
centered at $t$. In our approach, $h(t)$ measures the ``effective''
local density of the set $\Xi$ around the point $t$. We will discuss
this issue in the sequel. Right now, let us complete our basic
assumptions. We define the error kernel
\beqn
\label{errorfunction}
E(x,t):= \phi(x-t) - K(x,t),\quad x,t\in\Rd.
\eeqn
Notice that for each $t\in\Rd$, $E$ is a finite linear combination of
translates of $\phi$ using centers from $t\cup\Xi$.  We shall assume that

{\bf A4:} {\it There is a positive real
number $\nu>d$ and a constant $C>0$ depending only on $\phi$ such that
\beqn
\label{assumption1}
|E(x,t)|\le C h(t)^{\k-d} (1+{|x-t|\over
h(t)})^{-\nu},\quad x,t\in\Rd,
\eeqn
where $\k$ is the integer in {\bf A2}.}

As we will see in the examples that follow, the local density $h(t)$ must be chosen
to satisfy two properties: first, the ball $B_t(h(t)):=\{x\in\Rd:\ |x-t|\le
h(t)\}$ must contain a minimal number of centers from $\Xi$. This number is
determined by the parameter $\nu$. However, choosing $h(t)$ at this minimal
value leads to error kernel $E(\cdot,t)$ that are too large, forcing us to select
a large constant $C$. In such a case, it is usually preferable to increase
$h(t)$, so that more centers are captured in the enlarged ball $B_t(h(t))$.
By playing this game correctly at all points $t$, we can control the global
constant $C$. Needless to say, this comes at a price, since the density
function $h$ will enter our error bounds as well. 

%

\medskip
The remainder of this section will discuss two examples where the assumptions
{\bf A1-4} are satisfied.  These examples will provide a better understanding of
the assumptions as well as of
the nature of the smoothness spaces $W$ and the density function $h$.

\medskip\noindent
{\bf Example 1: Univariate splines.}\\  We consider the truncated power 
$\phi(t):=t_+^{\k-1}$ defined on
$\R$. We have the elementary and well-known representation
\beqn
\label{bspline}
f(x)=\frac{1}{(\k-1) !}\int_{-\infty}^\infty f^{(\k)}(t)\phi(x-t) \, dt,
\eeqn
which holds for all functions in $C_0^\k(\R)$. This means that {\bf
A2} holds for
$T:=D^\k/(\k-1)!$.  Now given $t\in\R$, let
$\Xi(t)=\{\xi_1(t),\dots,\xi_\k(t)\}\subset \Xi$ be the set   of the $\k$  
points in $\Xi$    that are closest to
$t$. The divided difference
\beqn
\label{Bspline}
[t,\xi_1(t),\dots,\xi_\k(t)]\phi(x-\cdot)=:
a(t)\phi(x-t)-\sum_{\xi\in\Xi(t)}A_0(t,\xi) \phi(x-\xi)
\eeqn
is the B-spline $x\mapsto M(x)$ associated with the knots
$\{t,\xi_1(t),\ldots,\xi_\k(t)\}$. Thus, we can take 
$A(t,\xi)=A_0(t,\xi)/a(t)$ for $\xi\in
\Xi(t)$ and $A(t,\xi)=0$ otherwise.   It follows that 
%
 $\displaystyle{h(t)=\max_{j=1,\dots,\k}
|t-\xi_j(t)|}$.
A well-known property of divided differences  (see \cite{DeLo}, p.
121) gives 
\beqn
\label{bcoeff}
|a(t)|= \prod_{j=1}^\k |t-\xi_j(t)|^{-1}\ge h(t)^{-\k}.
\eeqn
Since $|E(x,t)|= M(x)/|a(t)|\le h(t)^\k 
M(x)$, assumption
{\bf A4} follows from the facts that
$M(x)\le 1/h(t)$ and that
$M(x)$ vanishes for $x\notin (t-h(t),t+h(t))$.
As to the constant $C$, it can be chosen as $C:=4$ (for $\nu:=2$).
The verification of {\bf A3} is straightforward (via {\bf A1}).

It is worth stressing the fact that our ability to provide tight
error estimates for univariate splines is not only due to the banded
structure of the error kernel $E$. It was also due to the fact that
univariate spline theory tells us that $M(x)\le 1/h(t)$: thus,
while the actual coefficient $A(\cdot,\xi_i)$ of the truncated
power $(\cdot-\xi_i)^{\k-1}_+$ in the representation of $M(x)$ can
be arbitrarily large, the size of this coefficient does not affect
$\norm{M}_\infty$. Unfortunately, we are not aware of a multivariate
counterpart of this result. 
\endproof

\medskip\noindent
{\bf Example 2: Surface splines.} 
\\
The   multivariate analog of the truncated
power is the surface spline (see (\ref{surface}); it is
also known as the polyharmonic spline.) The best known surface spline is the 
bivariate {\it thin-plate spline}
$$\phi=|\cdot|^2\log|\cdot|,$$
which is, up to a constant, the fundamental solution for $\Delta^2$ 
when $d=2$.
               
We have noted earlier that property {\bf A2} holds for surface 
splines $\phi$ of any dimension.   We will analyse in detail
{\bf A4}, and will briefly discuss {\bf A3}. 

Let $\nu>d$ be the number that appears in 
{\bf A4}, and define
$$n:=\k-d+\nu.$$
Let $P$ be the space of all polynomials of degree $< n$ in
$d$ variables. Given a finite set $Z\subset\Rd$, we denote by $\Lam_Z$
the span of the functionals
$$\del_z\in P',\ z\in Z,\quad \del_z(p):=p(z), \ p\in  P.$$
Every $\lam:=\sum_{z\in Z}a(z) \del_z\in\Lam_Z$ extends to $C(\Rd)'$,
with the norm of the extension being
$$\norm{\lam}=\sum_{z\in Z}|a(z)|.$$

Now, for every $t\in\Rd$, we select a finite subset $\Xi(t)\subset\Xi$
that satisfies the following properties:

\begin{enumerate}

\item $h(t):=\max\{|t-\xi|:\ \xi\in \Xi(t)\}$ is ``as small as possible''.

\item There exists $\lam_t\in \Lam_{\Xi(t)}$ such that $\lam_t$ agrees with
$\del_t$ on $ P$.

\item $\norm{\lam_t}\le C$, for some $t$-independent constant $C$ that we 
choose in advance.

\end{enumerate}

Let us first remark that there are always sets $\Xi(t)$ satisfying
(ii) and (iii).%
\footnote{A sketch of the argument is as follows.
First, the claim is definitely correct when $\Xi=\Zd$. Therefore, there
exists $\del>0$, such that the claim is correct, provided that $\Xi$ has
non-zero intersection with any ball $B_a(\del)$, $a\in\Zd$.  Since
$\norm{\lambda}$ is invariant under dilation, the claim is thus correct
provided that $\Xi$ has non-zero intersection with any ball $B_a(R(1))$,
$a\in \frac{R(1)}{\del}\Zd$, and with $R(1)$ as in {\bf A1}. Thanks to 
assumption {\bf A1}, our $\Xi$ satisfies, indeed, this last property.
The argument as here is of mostly theoretical value, since
it employs a localization process that involves only a small subset of $\Xi$,
and results therefore in a density function that is prohibitively large.} %
We have said nothing about the size of the constant $C$ in (iii).
We do not  provide specific algorithms for 
choosing optimally $C$. The general rule of thumb is that by choosing $\Xi(t)$
to contain the $r+\dim P$ points in $\Xi$ that are closest to $t$,
with $r$ a small positive integer, we
should be able to find (with high probability, for a generic distribution
of centers) $\lam_t$ that satisfies the above.


Once we have chosen the functional 
$\lam_t=:\sum_{\xi\in\Xi(t)}A(t,\xi)\del_\xi$, we define
$$A(t,\xi)=0,\quad \hbox{on }\ \Xi\backslash\Xi(t),$$
and define the kernels
$$K(\cdot,t):=\sum_{\xi\in\Xi}A(t,\xi)\phi(\cdot-\xi),\	\quad
E(\cdot,t):=\phi(\cdot-t)-K(\cdot,t).$$

Next, recall that (up to a constant that depends only on $d$ and $\k$)
$\phi=|\cdot|^{\k-d}L$, with $L=\log|\cdot|$ whenever $\k-d$ is an
even integer, and $L=1$ otherwise.

     %
     %
      %
      
We now complete the proof of property {\bf A4}. Let $\Xi$ and $t$ be given.  
While we must verify {\bf A4} for every $t$ and every set $\Xi$ satisfying 
(i-iii), we can (by translating both $t$ and $\Xi$) assume without loss of 
generality that $t=0$.   Let $h:=h(t)$ be as in (i).  Suppose first that 
$|x|>2h$.  If $R$ is any polynomial of degree $<n$, then  
\beqn
\label{verify1}
|\phi(x)-K(x,0)|=|\phi(x)-\lam_0(\phi(x-\cdot))|\le
|\phi(x)-R(x)|+|\lam_0(\phi(x-\cdot)-R(x-\cdot))|.
\eeqn
In particular, choosing $R$ as the Taylor polynomial of degree $n-1$ at $x$ 
of $\phi$, the first term on the right side of \eref{verify1} is zero and we obtain
$$|\phi(x)-K(x,0)|\le C  \|\phi-R\|_{L_\infty(B)}$$
where $C$ is the constant in (iii) and $B$ is the ball of radius $h$ about $x$.  From the Taylor remainder formula we obtain
 \beqn
 \label{verify2}
 |\phi(x)-K(x,0)|\le  \frac{C|\phi|_{W^n(L_\infty(B))}h^n}{n !}\le
 CC'|x|^{\k-d-n}h^{n}\le CC''h^{\k-d}[1+\frac{|x|}{h}]^{-\nu},
 \eeqn
 where the constants $C',C''$ depend on $\phi$ and $\nu$  but are independent of $\Xi,t,x$.  Here, we used the fact that $|x|$ and $|x'|$ are comparable for $x'\in B$ because $|x|>2h$.
 
 If $|x|\le 2h$, then $|x-\xi|\le 3h$ for every $\xi\in\supp\lam_0$.
Assume momentarily that $h=1$. Then we simply estimate
\beqn
\label{verify3}
|\phi(x)-K(x,0)|=|\phi(x)-\lam_0(\phi(x-\cdot))|\le
(1+\norm{\lam_0})\norm{\phi}_{L_\infty(B)}\le CC''',
\eeqn
where now $B$ is the ball of radius $3$ about the origin, $C'''$ depends
only on $\phi$ and $n$, and $C$ is the constant in(iii).

Now, suppose that $h:=h(0)\not=1$.
Then,
dilating $\phi$, $\Xi$ and $\lam_0$ by  $h$, 
we note that
$$\phi(x/h)-\lam_0(\phi((x-\cdot)/h))=
h^{d-\k}(\phi(x)-\lam_0(\phi(x-\cdot)))+
(q(x)-\lam_0(q(x-\cdot))),$$
with $q$ a polynomial of degree $\le \k-d$ (viz., $q=0$ for odd $\k-d$, and
$q=-|\cdot|^{\k-d}h^{d-\k}\log h$ otherwise). 
Since $\delta_0-\lam_0$ annihilates
all such polynomials, we conclude that
\beqn
\phi(x/h)-\lam_0(\phi((x-\cdot)/h))=
h^{d-\k}(\phi(x)-\lam_0(\phi(x-\cdot))).\label{dil}\eeqn
Invoking now the analysis of the  ($h=1$)-case, we have that
$$|\phi(x/h)-\lam_0(\phi((x-\cdot)/h))|\le CC''',
$$
and hence, by (\ref{dil}),
$$|\phi(x)-\lam_0(\phi(x-\cdot))|\le CC'''  h^{\k-d}, 
$$
provided that $|x|\le 2h$. Altogether, for $x\in \Rd$,
 \beqn
 \label{verify5}
 |\phi(x)-K(x,0)| \le CC_1h^{\k-d}[1+\frac{|x|}{h}]^{-\nu},
 \eeqn
 with $C$ as in (iii), and $C_1$ a universal constant.

For general $t\in\Rd$, an argument identical to the above leads to
$$|\phi(x-t)-K(x,t)|\le C C_1 h(t)^{\k-d} 
(1+|x-t|/h(t))^{-\nu}.$$
This validates {\bf A4}. It also shows the constant that appears in {\bf A4}
is the  product of a constant  that is independent of $\Xi,x,t$ by
the uniform bound for the norms $\norm{\lam_t}$, $t\in\Rd$.

Concerning property {\bf A3}, for each fixed $\xi$, the function $A(\cdot,\xi)$ has compact support because of assumption {\bf A1} and our remarks above about the choice of $\Xi(t)$.  Since this function is also uniformly bounded as we have shown in the discussion of (iii), we see that {\bf A3} is also satisfied.
\endproof


\section{Approximation with a prescribed set $\Xi$ of centers}
\label{sect3}

In this section, we assume that the set $\Xi$ of centers is 
fixed in advance. We work under the 
assumption that $\Xi,\phi$ satisfy {\bf A1-A4}.   We shall prove a theorem for the  approximation of a 
given function $f\in W(L_p(\Rd),\phi)$   by the elements of
$S_\Xi(\phi)$.   In \S\ref{seclower}, we  will
extend the results of this section to   more general functions in $L_p(\Rd)$.

Since our goal is to derive error estimates that reflect the local density of
$\Xi$,  it may seem that we can employ our measure of density
$t\mapsto h(t)$ in such estimates.
However, it turns out that
$h$ may change too rapidly to allow effective error analysis (unless one 
replaces {\bf A4}  by the stronger assumption that $E(x,t)$ is 
supported in  the  domain $\{(x,t): |x-t|\le C h(t)\}$. However, the only
interesting example that satisfies this stronger condition is univariate
splines). To circumvent this difficulty, we introduce a companion density 
function $H$ that varies more slowly than the original $h$.

Given $\Xi, \phi$ and a local density $h$ that satisfy {\bf A1-A4}, 
we define
\beqn
\label{maximaldensity}
H(x):=\sup_{t\in\Rd} h(t)\left(1+\frac{|x-t|}{h(t)}\right)^{-r},\quad
x\in\Rd,
\eeqn
where $r$ is any fixed number satisfying
\beqn
\label{rcond}
0<r<\frac{\nu -d}{\k},
\eeqn
and with $\nu$ as in assumption (\ref{assumption1}).
The larger the value of $r$, the smaller the density function $H$ and the
better the estimates that we shall obtain.\footnote{It is therefore 
natural to try to take $r=\frac{\nu -d}{\k}$ in the analysis we give below;  
but this fails to work (barely).  One could introduce logarithmic factors in 
the definition of $H$ and get slightly improved results but at the expense of
notational complications that we want to avoid.}  Notice also that in the 
examples in the previous section, the number $\nu$ appearing in {\bf A4} can 
be chosen arbitrary large; however, the constant that appears in {\bf A4}
and the density function $h$ depend both on the selection of this $\nu$.

 We    assume   that $f\in W(L_p(\Rd),\phi)$ and then derive an
error estimate for approximating   $f$ by the elements of
$S_\Xi(\phi)$.  We first want to enlarge the space $S_\Xi(\phi)$ to include certain
infinite sums.  Given $1\le p\le \infty$, we define 
$$S_\Xi(\phi)_p$$ to be
the closure of $S_\Xi(\phi)$ in the topology of convergence in $L_p(\Omega)$
for each compact $\Omega\subset \Rd$.

 To describe the approximation procedure we are going to use, we recall the
kernel
$K$ (given by \eref{approxkernel}) which describes how $\phi(\cdot-t)$ is
approximated.  Also recall the error bound \eref{assumption1} for
$E:=\phi(\cdot-t)-K(\cdot,t)$  which we assume to hold for this
approximation.

  Given any positive weight function $w$ defined on $\Rd$,  we define the norm
  \beqn
  \label{wnorm}
  \|g\|_{L_p(w)}:=\|wg\|_{L_p(\R^d)},\quad 1\le p\le \infty ,
  \eeqn
  and the 
approximation error
\beqn
\label{approxerror}
\cE(f,S_\Xi(\phi))_{L_p(w) }:=\inf_{S\in
S_\Xi(\phi)_p}\|f-S\|_{L_p(w)},\quad 1\le p\le \infty .
\eeqn
Notice that $\|\cdot\|_{L_p(w)}$ differs from the more usual definition of 
weighted $L_p$-norms.

\begin{theorem}
\label{firsttheorem}
Suppose that $\Xi, \phi$ satisfy {\bf A1-4}, and let $1\le p \le \infty$.
If
$f\in W(L_p(\Rd),\phi)$ and
$w:=H^{-\k} $ , then 
 we have
\beqn
\label{firsterror}
\cE(f,S_\Xi(\phi)_p)_{L_p(w)}\le C_0 |f|_{W(L_p(\Rd),\phi)}
\eeqn
with $C_0=CC'$, $C'$ being dependent only on $\phi$, $\nu$ and $r$, and
$C$ is the constant that appears in {\bf A4}.
\end{theorem}
The parameter $\k$ (which was introduced in  Assumptions {\bf A2, A4})
determines the rate of decay of the error. Note that $\k$ appears
on both sides of (\ref{firsterror}) (in the definition of $L_p(w)$,
as well as in the definition of $W(L_p(\Rd),\phi)$). 

The theorem provides local error estimates in terms of the density
$H$.  Where $H$ is small, i.e. $\Xi$ is dense, the approximation bound is
better. The local nature of the error estimates is best captured
in the case $p=\infty$:

\begin{cor}
In the notations and assumptions of Theorem \ref{firsttheorem},
the error bound in the case $p=\infty$ can be restated as follows:
for every compact $\Ome\subset\Rd$ and every
$\eps>0$, there exists $g\in S_\Xi(\phi)$, such that,
for every $x\in \Ome$,
\beqn
\label{corfirsterror}
|f(x)-g(x)|\le \eps+ C_0 H(x)^\k |f|_{W(L_p(\Rd),\phi)}.
\eeqn
The constant $C_0$ is independent of $f,g,\eps,\Ome$ and $x$.
\end{cor}

{\bf Proof of Theorem \ref{firsttheorem}.} We begin by assuming that $f\in
C_0^\k(\Rd)$ and later use a completion argument to derive the general
case. We shall establish the estimate
\eref{firsterror} for $p=1,\infty$ and then derive the general case by
interpolation. For any function $g\in L_1(\Rd)+L_\infty(\Rd)$, we define 
\beqn
\label{operator1}
L(g,x):=\int_{\Rd} H(x)^{-\k} g(t)E(x,t)\, dt.
\eeqn
 Then $L$ is a linear operator and we shall show that it maps $L_p$
boundedly into itself for   $p=1,\infty$.  Once this is established, the 
Marcinkiewicz interpolation theorem implies that $L$ maps $L_p(\Rd)$ boundedly 
into itself for all
$1\le p\le \infty$.

First consider the case $p=1$.  We invoke the estimate \eref{assumption1}
and the definition of $H$  given in \eref{maximaldensity} to find
\beqn
\label{est1}
H(x)^{-\k}|E(x,t)|\le CH(x)^{-\k} h(t)^{\k-d} \left(1+{|x-t|\over
h(t)}\right)^{-\nu}\le Ch(t)^{-d}\left(1+{|x-t|\over
h(t)}\right)^{-\nu +r\k}.
\eeqn
Thus
\beqn
\label{est11}
\begin{array}{ll}
\|L\|_{L_1(\Rd)}&=\sup_{t\in \Rd}\int_{\Rd}H(x)^{-\k}|E(x,t)|\, dx\cr
&\le
\sup_{t\in \Rd} C
\int_{\Rd}h(t)^{-d} \left(1+{|x-t|\over h(t)}\right)^{-\nu+r\k}\,
dx  \\
&=C\int_{\Rd}(1+|y|)^{-\nu+r\k}\, dy\le C', 
\end{array}
\eeqn
where we used the fact that $-\nu+r\k<-d$.  

For the case $p=\infty$, we fix $x\in\Rd$ and define for each $j\in \Z$
the set
\beqn
\label{zj}
\Omega_j:=\{t\in \Rd: 2^{j-1}\le \frac{h(t)}{H(x)}< 2^j\}
\eeqn
Then,
\beqn
\label{est111}
 \int_{\Rd}H(x)^{-\k}|E(x,t)|\, dt=
\sum_{j\in\Z}\,\int_{\Omega_j}H(x)^{-\k}|E(x,t)|\, dt=:\sum_{j\in\Z} I_j.
\eeqn
We can estimate each of the integrals $I_j$ appearing in the sum by using
\eref{assumption1} to obtain
\beqn
\label{ij}
I_j\le
C2^{j\k}\int_{\Omega_j}[2^jH(x)]^{-d}\left(
1+\frac{|x-t|}{2^jH(x)}\right)^{-\nu}\,dt=C2^{j\k}
\int_{\Omega_j'}(1+|y|)^{-\nu}\, dy,
 \eeqn
where $\Omega_j'=[2^jH(x)]^{-1}(x-\Omega_j)$.
Since $\nu>r\k+d>d$, it is clear that $\sum_{j\le 1}I_j\le C'$.   For
$j> 1$, we use the definition of $H$ to find 
\beqn
\label{ij1}
\left(1+\frac{2|x-t|}{2^jH(x)}\right)^{r}\ge
\left(1+\frac{|x-t|}{h(t)}\right)^{r}\ge
\frac{h(t)}{H(x)}\ge 2^{j-1},\quad t\in \Omega_j.
\eeqn
In other words,
\beqn
\label{ij2}
\frac{|x-t|}{2^jH(x)}\ge  
\frac{2^{(j-1)/r}-1}{2} =:a_j,\quad t\in \Omega_j.
\eeqn
This means that 
\beqn
\label{ij3}
I_j\le C2^{j\k} \int_{|y|\ge a_j}(1+|y|)^{-\nu}\,dy\le C 2^{j\k}
a_j^{-\nu+d}
\eeqn
Since $\k<\frac{\nu-d}{r}$, we have that $\sum_{j>1}2^{j\k}a_j^{-\nu+d}$ is
finite.  This, together with (\ref{verify2}), yields
\beqn
\label{ij4}
\|L\|_{L_\infty(\Rd)} =\sup_{x\in\Rd} \int_{\Rd}H(x)^{-\k}|E(x,t)|\, dt\le
C.
\eeqn
Consequently, we have proved that $L$ boundedly maps $L_p$ into itself  for every
$1\le p\le \infty$.

\smallskip
Now, if  $f\in  W(L_p(\R^d),\phi)$, then by the definition of this space,
$Tf\in L_p(\R^d)$.  Hence, from what we have already proved,
\beqn
\label{proved1}
\|L(Tf)\|_{L_p(\R^d)}\le C_0\|Tf\|_{L_p(\R^d)}=C_0|f|_{W(L_p(\R^d),\phi)},
\eeqn
with $C_0=CC'$, with $C'$ an absolute constant, and $C$ the constant
that appears in {\bf A4}. 
Assume next that $f\in C_0^\k(\R^d)$ and define 
\beqn
\label{defF}
F:=\int_{\Rd}Tf(t)K(\cdot,t)\, dt.
\eeqn
  From {\bf A1-2} and the fact that $\supp Tf$ is compact, we deduce
  that the sum that defines $K(\cdot,\cdot)$ is finite on 
  $\Rd\times \supp Tf$. It follows then that 
\beqn
\label{repK}
F=\sum_{\xi\in\Xi}a(\xi)\phi(x-\xi), \quad
a(\xi):=\int_{\Rd}A(t,\xi)Tf(t)\, dt.
\eeqn
As said, the sum that defines $F$ is actually finite. Thus,
$F\in S_\Xi(\phi)$.

 Continuing under the assumption that $f\in C_0^\k(\Rd)$, we have
\beqn
\label{repK1}
H(x)^{-\k}[f(x)-F(x)] = H(x)^{-\k}\int_{\Rd} Tf(t)E(x,t)\, dt=L(Tf)(x).
\eeqn
It follows from what we have already proved that
\beqn
\label{repK2}
\cE(f,S_\Xi(\phi))_{w,p}\le \|f-F\|_{L_p(w)}\le C_0\|Tf\|_{L_p(\Rd)}=C_0|f|_{W(L_p(\Rd),\phi)}
,\quad 1\le p\le \infty.
\eeqn

\medskip
We now want next to extend \eref{repK2} to all of $W(L_p(\Rd),\phi)$.  Fix
$p\in [1,\infty]$, and let $f\in W(L_p(\Rd),\phi)$.    By the definition of
$W(L_p(\Rd),\phi)$ there is a sequence  of compactly supported functions
$f_n$, $n=1,2,\dots$,  from $C_0^\k$  such that  $f_n\to f$ in the norm of $W(L_p(\Rd),\phi)$.  Let $F_n$ be defined by \eref{defF}
for $f_n$.   We know that each of these $F_n$ is in $S_\Xi(\phi)$.  For any compact set $\Omega$, we have   $w(x)\ge c_\Omega>0$, $x\in\Omega$,  and therefore by writing $F_m-F_n=F_m-f_m-(F_n-f_n)+(f_m-f_n)$ we find
\begin{eqnarray}
\label{Cauchy}
\|F_m-F_n\|_{L_p(\Omega)}&\le& \|f_m-f_n\|_{L_p(\Omega)}+c_\Omega^{-1}\|f_m-f_n-(F_m-F_n)\|_{L_p(w)}\cr
&\le& \|f_m-f_n\|_{L_p(\Rd)}+C_\Omega\|T(f_m-f_n)\|_{L_p(\Rd)} ,
\end{eqnarray}
where the last inequality uses \eref{repK2} for $f=f_m-f_n$.  This shows that $(F_n)$
is a Cauchy sequence in the topology of $L_p$-convergence on compact sets.  By definition, its limit $G$ is in $S_\Xi(\phi)_p$.  Again, for any compact set $\Omega$ in $\Rd$,  we have
\beqn
\label{finally}
\|w(f-G)\|_{L_p(\Omega)}\le \lim_{n\to \infty}\|f_n-F_n\|_{L_p(w)}\le C_0\lim_{n\to \infty}\|T(f_n)\|_{L_p(\Rd)}=C_0\|T(f)\|_{L_p(\Rd)},
\eeqn
with $C_0$ the constant of \eref{repK2}.
Since $\Omega$ is arbitrary, we find
\beqn
\label{finally1}
\|w(f-G)\|_{L_p(\Rd)}\le C_0\|T(f)\|_{L_p(\Rd)}.
\eeqn
Since $G\in S_\Xi(\phi)_p$, we can replace the left side of \eref{finally1} 
by $\cE(f,S_\Xi(\phi)_p)_{L_p(w)}$.  This completes the proof of the theorem.
\endproof


\bigskip
Theorem \ref{firsttheorem} deals with the approximation of functions that are
optimally smooth, i.e., in the space $W(L_p(\Rd),\phi)$. In
\S\ref{seclower}, we establish results concerning the approximation of
functions that are less smooth. For such functions, the weight
$H^{-\k}$ is too strong. To this end, we state a counterpart of Theorem
\ref{firsttheorem} for mollified versions of the original weight
$w$. We still assume here that $f$ is optimally smooth. This assumption
will be dropped in \S\ref{seclower}.

Suppose that $0< s< \k$.  We continue to work under the assumptions {\bf A1-4}.   If $f\in W(L_p(\Rd),\phi)$ then
\beqn
\label{Wsseminorm}
\|h^{\k-s}Tf\|_{L_p(\Rd)}
\eeqn
 is finite because $h$ is bounded (see {\bf A1 (ii)}).
 
\begin{theorem}
\label{secondtheorem}
If  $0<s< \k$ and  $f\in  W(L_p(\Rd),\phi)$ (for some $1\le p\le \infty$),
then for $w:=H^{-s}$ we have
\beqn
\label{seconderror}
\cE(f,S_\Xi(\phi))_{L_p(w)}\le C_0 \|h^{\k-s}Tf\|_{L_p(\Rd)}.
\eeqn
with $C_0$ as in Theorem \ref{firsttheorem}.
\end{theorem}

{\bf Proof:}  The proof is very similar to the proof of Theorem
\ref{firsttheorem}.  We  first remark that  the two bounds
\beqn
\label{se1}
\sup_{t\in \Rd}\int_{\Rd}H(x)^{-s}|E(x,t)|h(t)^{s-\k}\, dx\le C,\quad 
\sup_{x\in \Rd}\int_{\Rd}H(x)^{-s}|E(x,t)|h(t)^{s-\k}\, dt\le C,
\eeqn
hold with $C$ a constant depending only on $d,s$.  Indeed, thanks to Assumption
{\bf A4} it is sufficient to prove the above boundedness with the integrand
replaced by 
$$H(x)^{-s}h(t)^{\k-d}h(t)^{s-\k} \left(1+{|x-t|\over
h(t)}\right)^{-\nu}= H(x)^{-s}h(t)^{s-d} \left(1+{|x-t|\over
h(t)}\right)^{-\nu}.$$
Since we assume  $\nu>\k r+d\ge sr+d$, the argument as given in the proof of  
Theorem \ref{firsttheorem} applies here {\it verbatim} to yield
(\ref{se1}) .

The bounds in (\ref{se1}) now imply that the linear operator
$$
L(g):= H(x)^{-s}\int_{\Rd}g(t)E(x,t)h^{s-\k}(t)\, dt$$
is bounded on $L_p(\Rd)$ for $p=1,\infty$.  By interpolation, we derive
that this operator is bounded on $L_p(\Rd)$ for all $1\le p\le \infty$.
Using this for $g=h^{\k-s}T(f)$, we derive \eref{seconderror} for all $f\in
W(\phi,L_p(\Rd))$ in the same way we have proven Theorem \ref{firsttheorem}.  
\endproof

\section{Wavelet decompositions}
\label{secwave}
In the remaining sections of this paper, we shall be in need of a local
multiscale basis on $\R^d$.  We shall employ a standard multivariate wavelet
basis for this purpose. This basis will be used only as a tool for proving
various results.  In this section, we recall the form of such a basis and
some of its properties which will be important to us. In particular, we shall
need its characterization of Triebel-Lizorkin spaces.  There are
several books that discuss wavelet decompositions and their characterization
of these spaces (see e.g. \cite{Meyer}).  We 
also refer to the article of Daubechies \cite{D92} for the construction of 
wavelet bases of the type we want to use. 


Let $\cD$ denote the set of dyadic cubes in $\Rd$ and $\cD_j$ the set of
dyadic cubes of side length $2^{-j}$ (thus
$\cD=\cup_{j=-\infty}^\infty\cD_j$).  Each
$I\in\cD_j$ is of the form
\beqn
\label{dyadic}
I=2^{-j}[k_1,k_1+1]\times\cdots \times 2^{-j}[k_d,k_d+1]=2^{-j}(k+[0,1]^d),\quad
k=(k_1,\dots,k_d)\in\Zd.
\eeqn
For each such $I$, we denote its side length by $\ell(I)$:
$$\ell(I):=2^{-j},\quad\forall I\in\calD_j.$$
Finally, let
$$E:=\{1,\ldots ,2^d-1\},\quad \calV:=\calD\times E.$$
Given
$\v=(\iv,e_\v)$ ($\iv\in \calD$, $e_v\in E$), we denote 
$$\ell(\v):=\ell(I_\v),$$
and by
$$|v|:=\ell(v)^d$$
the volume of the cube $\iv$.

A {\it wavelet basis} is an orthonormal basis for $L_2(\Rd)$ with particular
structure and properties. The wavelets are indexed by the set
$\calV$:
$$w_\v,\quad \v\in\calV.$$

Each wavelet $w_{\v}$, with $\iv=2^{j}(k+[0,1]^d)\in\calD$, is supported in a cube 
$\bar \iv$, with
$${\bar \iv}=2^{j}(k+A_0[0,1]^d),$$
with $A_0$ some fixed constant that depends only on the specifics of
the wavelet system we choose.

We normalized initially the wavelet system in $L_2(\Rd)$. The $L_p$-norm
of the $p$-normalized wavelets, $1\le p\le \infty$:
\beqn
\label{waveletf}
\psi_{\v,p}:=|\v|^{{1\over 2}-{1\over p}}w_\v
\eeqn
depends only on their {\it type}, i.e., on the index $e_\v\in E$. 
%
%
Each locally integrable
function $f$ defined on $\Rd$ has a wavelet decomposition
\beqn
\label{wavedec}  
 f=\sum_{\v\in\calV}f_\v\psi_{\v,p},\quad f_\v:=f_{\v,p}:=\langle
f,\psi_{\v,p'}\rangle,\quad 1/p+1/p'=1.
\eeqn
Here $f_{\v}$ depends on the $p$-normalization that has been
chosen  but the product $f_\v\psi_{\v,p}$ is independent of $p$.   In this
paper, it will be convenient to normalize the wavelets in $L_\infty$:
$$\psi_\v:=\psi_{\v,\infty}=|\v|^{{1\over 2}}w_\v.$$
The series \eref{wavedec} converges 
absolutely to $f$ in the $L_p$-norm in the  case
$f\in L_p(\Rd)$ and 
$1\le p<\infty$, with $H_1(\Rd)$ replacing $L_1(\Rd)$,
and conditionally  in the case $p=\infty$ with
$L_\infty(\Rd)$ replaced by $C(\Rd)$.

%
%
One of the most important properties of wavelet systems (and the one we
need in this paper) is the characterization of smoothness spaces
in terms of the wavelet decomposition. This means that we can use 
the wavelet decomposition in order to {\it define} those smoothness
spaces.  For the definition of  Triebel-Lizorkin spaces $F_{p,q}^s$
in terms of 
wavelet coefficients, we fix $s>0$.
Then, for
$0<p,q<
\infty$,  
 \beqn
 \label{TLnorms}
 |f|_{F_{p,q}^s}:= \| M_{q}(f)\|_{L_p(\Rd)},
 \eeqn
 where
 \beqn
 \label{maximal}
M_s(f)(x):=M_{s,q}(f)(x):=\left(
 \sum_{\v\in\calV} \ell(\v)^{-qs}|f_{\v,\infty}|^q\chi_{\bar I_\v}(x)
 \right)^{1/q}.
 \eeqn
The definition does not depend on the wavelet system we choose, provided
that the wavelets are $m$-times differentiable, and 
have $m$ vanishing moments (viz., their Fourier transform has a $m$-fold
zero at the origin), for sufficiently large $m$.\footnote{The basic
requirement is that $m>s$. Additional requirements are imposed
in case $p<1$ or $q<1$. For us, the only thing that matters is the existence
of some wavelet system that can be used to define the Triebel-Lizorkin space
$F^s_{p,q}$. In particular, we can, and do, allow the wavelet system
to depend on $s,p,q$.}

We should also make some specific
remarks about our definition. In our  definition we have defined the maximal
function with $\chi_{\bar \iv}$ where $\bar \iv$ the support cube of the 
$\psi_\v$.  The usual definition of $F_{p,q}^s$  uses $\chi_{\iv}$ 
instead of $\chi_{\bar \iv}$.  It is easy to see that these two definitions 
give equivalent norms by using the Fefferman-Stein inequality 
mentioned in the proof of Lemma \ref{maximalfun}.  We also remark that 
Triebel-Lizorkin spaces are usually defined using Littlewood-Paley 
decompositions. Our definition agrees with the classical definition if
(the wavelet system is chosen appropriately,  cf.\ the above footnote)and
the space $F_{p,q}^s$ continuously embeds into $L_1$. 
       
The definition extends naturally to the $q=\infty$ case, with
$M_\infty(f)$ defined by
$$M_\infty(f)(x):= M_{s,\infty}(f)(x):=
\sup_{\v\in\calV} \ell(\v)^{-s}|f_{\v,\infty}| \chi_{\bar I_\v}(x).$$
While we followed so far the tradition of of assuming $p<\infty$,
we will use the space $F_{\infty,\infty}^s$, whose semi-norm is
$$|f|_{F^s_{\infty,\infty}}:=
\sup_{\v\in\calV} \ell(\v)^{-s}|f_{\v,\infty}| \chi_{\bar I_\v}(x),$$
and which usually appears in the literature as the Besov space
$B_{\infty,\infty}^s$.  Notice that $B_{\infty,\infty}^s$ is compactly embedded in $C(\R)$.
  
%
The
quasi-norm in
$F_{p,q}^s(\Rd)$ is defined by
\beqn
\label{TLnorm}
\|f\|_{F_{p,q}^s(\Rd)}:=\|f\|_{L_p(\Rd)}+|f|_{F_{p,q}^s(\Rd)}.
\eeqn

\section{Approximation of functions with lower smoothness}
\label{seclower}
The estimate \eref{seconderror} is unsatisfactory because it can be applied
only to a small subset of functions in $L_p$.  In this section, we shall 
remove this deficiency.  Our method for doing this is an \lq interpolation 
of operators' type argument which decomposes a   general function into a 
smooth part to which \eref{seconderror} can be applied and   a second 
nonsmooth part which is small.  We shall restrict our discussion to the 
case where the operator $T$ is an elliptic differentiable operator of order 
$\k$ with constant coefficients and $\phi$ is its fundamental solution (on 
$\Rd$). For such $T$, we can choose the wavelet system in the previous
section to satisfy
\beqn
\label{tpsi}
T(\psi_\v)(x)\le C'\ell(\v)^{-\k}\chi_{{\bar \iv}}(x),\quad \v\in \calV.
\eeqn

We fix $1\le p\le \infty$  and continue to work under assumptions {\bf A1-4}.
We recall in the present case 
$W(L_p(\Rd),\phi))$ is the Sobolev space $W^{\k}(L_p(\Rd))$.
Then, each  $f\in L_p(\Rd)$ has a series representation 
$$f=\sum_{\v\in\calV}f_\v\psi_\v,$$
with the sum convergent in $L_p(\Rd)$.    We shall work in this section 
exclusively with $L_\infty$-normalized wavelets $\psi_\v$.

Now, let $f$ be a function in the Triebel-Lizorkin space $F^s_{p,\infty} $, 
$s< \k$, and let $h$ be the density function for $\Xi$. We decompose 
$f=f_h^++f_h^-$ in the following way:
$$f_h^+:=\sum_{  \ell(\v)\ge h(\v)}\psi_\v f_\v,\quad f_h^-:=f-f_h^+,$$
where
$$h(\v):=\norm{h}_{L_\infty(\bar \iv)}.$$
 
We shall first estimate how well $f_h^+$ can be approximated by elements 
from $\cS_\Xi(\phi)$.
\begin{lemma}
\label{pluslemma}
Let  $1\le p\le \infty$ and  $0<s<\k$.  If $f\in F_{p,\infty}^s$, then for 
$w:=H^{-s}$, we have
\beqn
\label{plus1}
\dist(f_h^+,\cS_\Xi(\phi))_{L_p(w)}\le C(s,d)\|f\|_{F_{p,\infty}^s},
\eeqn
with a constant $C$ as in Theorem \ref{firsttheorem}.
\end{lemma}
{\bf Proof:} 
>From (\ref{tpsi}), we obtain for any $x\in\Rd$,
\begin{eqnarray}
|h^{\k-s}(x)T(f_h^+)(x)|&\le& h^{\k-s}(x)
\sum_{\ell(\v)\ge h(\v)}|f_{\v}| |T(\psi_{\v})(x)|\cr
&\le& C' h^{\k-s}(x)\sum_{\ell(\v)\ge 
h(\v)}\ell(\v)^{-\k+s}\ell(\v)^{-s}|f_{\v}|\chi_{\bar \iv}(x)\cr
&\le& 
C' M_{s,\infty}(f)(x)\, h(x)^{\k-s}
\sum_{\ell(\v)\ge h(x)} \ell(\v)^{-\k+s}\chi_{\bar \iv}(x) \cr
 &\le &
CM_{s,\infty}(f)(x) h(x)^{\k-s}  h(x)^{-\k+s}\le CM_{s,\infty}(f)(x).
\end{eqnarray}
Here we have used the fact that there is an absolute constant 
$C_1$ depending only on the support size of the wavelet such that for any 
dyadic level $j$, there are at most $C_1$ $I\in\cD_j$ which contain the 
given point $x$.  This means that the above series can be 
compared with a geometric series and can be bounded by a fixed multiple of 
its largest term. 
Thus, the constant  $C$   depends only  $s$ and $d$.
Therefore, by Theorem \ref{secondtheorem},
$$\dist(f_h^+,S_\Xi(\phi))_{L_p(H^{-s})}\le C
\norm{h^{\k-s}T(f_h^+)}_{L_p(\Rd)}\le C\norm{M_{s,\infty}(f)}_{L_p}=
C\norm{f}_{F^s_{p,\infty} },$$
which completes the proof of the lemma.
\endproof

\medskip{\bf Remark.} Formally, the proof given in the lemma does not
cover the case $p=1$. The reason is that the wavelet
representation of $f\in L_1(\Rd)$ does not always converge
to $f$. However, the lemma does extend to all $f\in L_1(\Rd)$,
provided that we use an inhomogeneous wavelet representation.
Such representation takes the form
$$f=\sum_{\ell(\v)<2^J}\inpro{f,\psi_\v}\,
\psi_\v+\sum_{\ell(\v)=2^J}\inpro{f,\tilpsi_\v}\tilpsi_\v=:f_1+f_2,$$
with the modified wavelets $\tilpsi_\v$ supported in exactly the same
cube $\biv$ as their original wavelet counterparts, and satisfy
$T(\tilpsi_\v)\le C'2^{-J\k}\chi_{\biv}$. The integer $J$ can be chosen
at will. This modified expansion converges for every $f\in L_1(\Rd)$.
We can use the above
inhomogeneous wavelet expansion in the proof of the Lemma,  
since we know (see the discussion on surface splines in \S2) that
the density function is {\it bounded}, which implies that
the term $f_h^+$ in the decomposition of $f$ contains
the entire expansion of the above $f_2$ (for a suitable large $J$ that
depends on the bound we have on $h$, but on  nothing else.) The argument
in the proof of the lemma can be then repeated  {\it verbatim} for the
case $p=1$. However, the smoothness space that is characterized by
the inhomogeneous expansion is the inhomogeneous Triebel-Lizorkin
space. For this reason, we have stated the result with respect to the
full norm $\|f\|_{F^s_{p,\infty}}$. For $1<p \le \infty$, the result is also
valid with $\|f\|_{F^s_{p,\infty}}$ replaced by
$|f|_{F^s_{p,\infty}}$.

\medskip
We are left with bounding the $L_p(H^{-s})$-norm of $f_h^-$.   
\begin{lemma}
\label{minuslemma}
Fix $1\le p\le \infty$. If $0<s<\k$ and $w:=H^{-s}$, then 
\beqn
\label{minus1}
\|f_h^-\|_{L_p(w)}\le C(s,d)\|f\|_{F_{p,\infty}^s},
\eeqn
where $C(s,d)$ depends only on $s$ and $d$.
\end{lemma}

{\bf Proof:}  Since $f_h^-=\sum_{\ell(\v)<h(\v)}f_\v\psi_{\v}(x)$
and $|\psi_\v(x)|\le C'\chi_{\bar \iv}(x)$ for an absolute constant $C'$, we have
\beqn
\label{minus2}|f_h^-(x)|\le CM_{s,\infty}(f)(x) \sum_{\ell(\v)<h(\v)} 
\ell(\v)^{s}\chi_{\bar  \iv}(x)\eeqn
  Given an $x\in \R^d$ for which $x\in \bar \iv$ and $ \ell(\v)<h(\v)$, we know 
  that there is a $t\in \bar \iv$ such that  $h(\v)=h(t)$.  For this $t$, we have 
 \beqn
 \label{minus3}
 H(x)\ge h(\v)(1+{|x-t|\over h(\v)})^{-r}.
 \eeqn
Since, $|x-t|\le \sqrt{d} \ell({\bar \iv})\le \sqrt{d} A_0\ell(\v)\le
\sqrt{d} A_0 h(\v)$,  we
find   that $h(\v)\le C H(x)$ with $C$ depending only on $r$.
Thus, 
\beqn
\label{minus4}
\sum_{\ell(\v)<h(\v)} \ell(\v)^{s}\chi_{\bar  \iv}(x)
\le 
\sum_{  \bar \iv\ni x,\,\ell(\v)\le CH(x)} \ell(\v)^{s}\le  
C'H(x)^{s},
\eeqn
where $C'$ depends on $s$ and $d$.
Here, as in the previous lemma we compared the series in \eref{minus4}  with a geometric series and  bounded it by a fixed multiple of 
its largest term. 
>From \eref{minus4}, we obtain
$$\norm{H^{-s}f_h^-}_{L_p}\le
C\norm{M_{s,\infty}(f)}_{L_p}=C\norm{f}_{F^s_{p,\infty}(\Rd)},$$
which proves the lemma.\endproof

Combining the two lemmas, we arrive at the following result.

\begin{theorem}
\label{wavetheorem}
 Assume that $\phi$ and $\Xi$ satisfy the assumptions {\bf A1-4} with
 respect to a homogeneous differential operator $T$ with constant
 coefficients (and degree $\k$). For every $1\le p\le \infty$, for every
$0<s< \k$, and for every $f\in F^s_{p,\infty} $ we have
\beqn
\label{wave1}
 \dist(f,S_\Xi(\phi))_{L_p(H^{-s})}\le C(s,d)\norm{f}_{F^s_{p,\infty} },
 \eeqn
 where $C$ as in Theorem \ref{firsttheorem}, 
and $H$ is the majorant density function associated with $\Xi$.
\end{theorem}

%
We can also derive results for approximation in Sobolev spaces, as well as
in Besov spaces. For the Sobolev space $W_p^s$, $1\le p<\infty$,
$0<s<\k$, we have that $W_p^s=F_{p,2}^s\subset F_{p,\infty}^s$,
\cite{Triebel}, and hence
Theorem \ref{wavetheorem} implies that
\beqn
\label{sobres}
 \dist(f,S_\Xi(\phi))_{L_p(H^{-s})}\le C\norm{f}_{W_p^s(\Rd)}. 
 \eeqn
As to Besov   spaces, 
since $B_p^s(L_p^{\vphantom s}(\Rd))$ is continuously embedded 
into $F_{p,\infty}^s$ (see \cite {Triebel}), we have, by the same theorem,
that, for $1\le p<\infty$
and $0<s<\k$,
\beqn
\label{wave11}
 \dist(f,S_\Xi(\phi))_{L_p(H^{-s})}\le C\norm{f}_{B_p^s(L_p^{\vphantom s}(\Rd))}. 
 \eeqn
 %
This latter statement extends trivially to $p=\infty$, since
$F_{\infty,\infty}^s=B_\infty^s(L_\infty)$.

\medskip
Note that we have restricted out attention to a differential operator
$T$. The reason is that our analysis depends on two properties of the
wavelet system. The first is the characterization of 
Triebel-Lizorkin spaces in terms of wavelet decompositions, and the second
is (\ref{tpsi}). In order to extend the result of this section to
more general $T$, we will need a representation system that will
satisfy, first and foremost, a bound analogous to (\ref{tpsi}) with respect 
to the more general $T$. We can then define using that system smoothness
spaces that are analogous to the Triebel-Lizorkin spaces
from  \S\ref{secwave}, and establish an analog of Theorem \ref{wavetheorem}.

\section{Nonlinear approximation}
\label{secnonlinear} 

 We shall now turn to a  
different setting.  We assume  $f\in L_p(\Rd)$ is a function 
that we wish to approximate using the shifts of $\phi$.  In contrast to
the problems studied so far where the set of centers is prescribed in advance, we shall now allow the choice of the centers  to be made dependent on $f$. We are interested in how well
we can approximate $f$ using at most $N$ such centers.

Let $\Sigma_N:=\Sigma_N(\phi)$   be the set of all functions $S$
for which there 
is a set $\Xi$ of cardinality $N$ such that
\beqn
\label{defSigma}
S=\sum_{\xi\in \Xi}a_\xi\phi(\cdot-\xi).
\eeqn
We then define the approximation error
\beqn
\label{nl}
\sigma_N(f)_p:=\inf_{S\in \Sigma_N} \|f-S\|_{L_p(\Rd)}.
\eeqn
This form of nonlinear approximation is known as {\it $N$-term
approximation}.    

Our setting is different than that considered in previous sections.
We begin with a function $\phi$ that satisfies assumption {\bf A2}.  For 
simplicity, we also suppose that the operator $T$ is a homogeneous 
differential operator of order $\k$.   The reader can easily abstract the conditions we use about $T$ to get a more general theorem.

We do not assume {\bf A1, A3, A4} since $\Xi$ is not specified in advance.
Rather, we shall assume that for any dyadic cube   there is a collection of
points near this cube that locally satisfy {\bf A3-4}.  To make this assumption
precise, we fix the order $\k$ of the differential operator $T$ and fix 
a wavelet system $\{\psi_{\v}\}_{\v\in\calV}$ of $\k$-times 
differentiable compactly supported wavelets of the form described in \S 
\ref{secwave}.  As in that section, for each $\v$, we denote by $\biv$ the 
smallest cube which contains the support of $\psi_\v$.
Recall that we know that $\biv$ has size comparable to that of $\iv$, 
namely $\ell({\biv})\le A_0\ell(\v)$  for a fixed 
constant $A_0$ that depends only on $\k$.

We make the following assumption (about $\phi$) in this section.

\medskip
{\bf A5:} {\it  Given $\nu>d$, there is an absolute constant $C_0$ and an 
integer $  N_0$ such that for any $N\ge N_0$, and any $\v\in \calV$,  there is 
a set $\Xi_{\v,N}\subset \Rd$ consisting of $N$ points with the 
following property:
%
  %
There is a linear combination     
 $$K_{\v,N}(\cdot,t)=\sum_{\xi\in \Xi_{\v,N}} A_{\v,N}(t,\xi)\phi(\cdot-\xi),
 \quad t\in {\biv},$$
    with  $A_{\v,N}(\cdot, \xi)$ in $L_\infty(\biv)$ such that 
    $\phi(x-t)-K_{\v,N}(x,t)$ satisfies 
\beqn
\label{apppsi22}
|\phi(x-t)-K_{\v,N}(x,t)|\le C_0 h_{\v,N}^{\k-d}
\left(1+\frac{|x-t|}{h_{\v,N}}\right)^{-\nu}, \quad t\in \biv,\
x\in\Rd,
\eeqn
}
with
$$
h_{\v,n}:=\frac{\ell(v)}{N^{1/d}}.
$$

\medskip
The new assumption {\bf A5} can be easily shown to be satisfied for surface
splines (Example 2 in \S\ref{sect2}). Indeed, given $N$, we choose
$m$ such that $m^d\le N$, and choose $\Xi_{\v,N}$ to be the vertices
of a uniform grid on $\bar \iv$ with mesh size $\tilh_{\v,N}:=\ell({\bar 
\iv})/(m-1)$.  Assuming that $m\ge \k-d+\nu$ (see Example 2),
it is then easy to follow the dilation argument given in 
Example 2 to conclude that the linear functionals $\lam_t$, $t\in {\bar \iv}$
are uniformly bounded, independently of ${\bar \iv}$ and $m$ (hence of
$N$). Also, $\tilh_{\v,N}$ as above satisfies $\tilh_{\v,N}\le C h_{\v,N}$.
Thus, (\ref{apppsi22}) follows from
the validity of {\bf A4} for surface splines (see Example 2 in \S
\ref{sect2}).  The constant $C$ depends
only on $d$ and the number $A_0$, hence is universal for all of our
purposes. Thus, $N_0$ in this case is $(\k-d+\nu)^d$. 
 
%
%
%

 To prove theorems about $N$-term approximation of a function $f$ by  
 shifts of the function $\phi$, we will  first   represent $f$ in a wavelet decomposition and approximate the individual terms in this decomposition.
We begin by seeing how well we can approximate the 
wavelet  $\psi_\v$ (normalized in $L_\infty$),
by using a budget of $N_{\v}$ centers.   Given $\v$, and any integer $N_{\v}\ge
N_0$, we let
$\Xi_{\v,N_{\v}}$ be the  set of points satisfying {\bf A5}. 
 Let us denote  
\beqn
\label{apppsi}
S_{\v,N_{\v}}:=\int_{\Rd}T(\psi_\v)(t)K(\cdot,t)\,dt,\quad
K(x,t):=\sum_{\xi\in\Xi_{\v,N_{\v}}}A(t,\xi)\phi(x-\xi).
\eeqn
Notice that since $T(\psi_\v)=0$, $t\notin \bar \iv$, the integral in
\eref{apppsi}  only goes over $t\in \bar \iv$ and therefore by  {\bf A5}, this integral is finite.

 Now, an important property of the wavelets that we need here is that
$\psi_\v$ satisfies
$\|T\psi_\v\|_{L_\infty}\le C_0\ell(\v)^{-\k}$. Also, thank to {\bf A2}, 
$\psi_\v=\int_{\Rd}T(\psi_\v)(t)\phi(\cdot-t)\, dt$. We 
use this together with {\bf A5} to derive  the following bound 

\begin{eqnarray}
\label{psiappr}
|\psi_\v(x)-S_{\v,N_{\v}}(x)|&\le&
C_0C_1
\left(
\frac{h_{\v,N_{\v}}}{\ell(\v)}\right)^
{\k}\int_{{\bar \iv}} 
\left(1+\frac{|x-t|}{h_{\v,N_\v}}\right)^{-\nu}h_{\v,N_{\v}}^{-d}\, dt\cr
&\le & CN_{\v}^{-\k/d}\left(1+\frac{dist(x, \bar 
\iv)}{\ell(\v)}\right)^{-\nu+d}.
\end{eqnarray}

Let us describe now our approximation algorithm and analyze its performance.
In the algorithm, we are given a budget $N$ of centers, and
invest a nominal amount $c_\v\ge 0$ in each
wavelet $\psi_\v$, $\v\in\calV$; we refer to $c_\v$ as {\bf 
cost}.  We ensure that the total cost $\sum_\v c_\v$ does not exceed the given
budget $N$. 
Since $c_\v$ may not be an integer, and since the minimal number of centers
that we can use is $N_0$, the cost $c_\v$ allows us to approximate
the term $f_\v\psi_\v$ in the wavelet expansion (\ref{wavedec})
of $f$ by investing
\beqn
\label{defni}
N_\v:=\lfloor c_\v\rfloor
\eeqn
centers, provided $c_\v\ge N_0$. In this case, (\ref{psiappr}) gives us the
estimate
\beqn
\label{tlIerror}
|\psi_\v(x)- S_{\v,N_{\v}}(x)|\le
CN_{\v}^{-\k/d}\left(1+\frac{dist(x, \bar \iv)}{\ell(\v)}\right)^{-\nu+d}
\le 
C'c_{\v}^{-\k/d}\left(1+\frac{dist(x, \bar \iv)}{\ell(\v)}\right)^{-\nu+d}.
\eeqn
If $c_\v<N_0$, we do not approximate $\psi_\v$ at all, and get then
\beqn
\label{tlIerror1}
|\psi_\v(x)- S_{\v,N_{\v}}(x)|
=|\psi_\v(x)|\le C  \chi_{\bar\iv}(x),\quad x\in\Rd,
\eeqn
since   the wavelets are uniformly
bounded.

Now, suppose that we are given a budget of $N$ centers, determine a cost
distribution $(c_\v)_\v$ (with $\sum_\v c_\v \le N$), and would like to estimate
the $L_p$-norm of the error when approximating the term 
$f_\v\psi_\v$ in the wavelet expansion of $f$ using the designated cost
$c_\v$. According to (\ref{tlIerror}) and (\ref{tlIerror1}), 
our error will be determined, up to
a universal constant, by the $p$-norm of
\beqn
\label{defR}
\sum_\v |f_\v| |\psi_\v-S_{\v,N_\v}|=:
\sum_\v |f_\v| R_\v.
\eeqn
Here,
$(f_\v)$ are the wavelet coefficients of the approximand $f$.

\begin{lemma}
\label{maximalfun}
Let $1\le p\le \infty$ and suppose that the constant $\nu$ appearing in {\bf A5} is $>2d$.  If $R_\v$ is defined as in {\rm \eref{defR}}, then
\beqn
\label{mf1}
\norm{\sum_{\v\in\calV}|f_\v|R_\v}_{L_p(\Rd)}\le
C\norm{\sum_{\v\in\calV} [\max (1,c_\v)]^{-\k/d}|f_\v|\chi_{\iv}}_{L_p(\Rd)}.
\eeqn
\end{lemma}
{\bf Proof:} This inequality can be derived using the method of proof for Fefferman-Stein inequalities \cite{FS}.  In fact, for $1\le p<\infty$, it can be derived   directly from these inequalities  as follows.  Let    $M_0$ be the Hardy-Littlewood maximal operator
\beqn
\label{FS1}
M_0(g)(x):=\sup_{Q\ni x}\frac{1}{|Q|}\int_Q|g(u)|\, du
\eeqn
where the supremum is taken over all cubes $Q$ that contain $x$.
Then, for $1<p<\infty$ and any real numbers $(a_\v)$,  the   Fefferman-Stein inequality says 
\beqn
\label{FS2}
\norm{\sum_{\v\in\calV}|a_\v |M_0(\chi_{ I_\v})}_{L_p(\Rd)}\le
C\norm{\sum_{\v\in\calV} |a_\v| \chi_{ I_\v}}_{L_p(\Rd)},
\eeqn
where $C$ depends only on $p$ as $p$ gets close to $1$ and $\infty$.
Now a direct calculation shows that for a constant $C_0$ depending only on $d$ we have
\beqn
\label{FS3}
M_0(\chi_{ I_\v})(x)\ge C_0\left(1+\frac{dist(x,  \iv)}{\ell(\v)}\right)^{-d} \ge C'_0\left(1+\frac{dist(x, \bar \iv)}{\ell(\v)}\right)^{-\nu+d} 
\eeqn
where, in the last inequality,  we used our assumption  that $\nu>2d$ and the fact that $1+\frac{\dist(x,I_\v)}{\ell(\v)}\le C(1+\frac{\dist(x,\bar I_v)}{\ell(\v)})$ for all $x\in \Rd$ with a constant $C$ depending only on the space dimension $d$. It follows from this and   \eref{psiappr}, \eref{tlIerror1}  that  
$$
R_\v(x)\le  C [\max(1,c_{\v})]^{-\k/d} M_0(\chi_{I_\v})(x),\quad x\in\Rd,
$$
with $C$ again depending only on $d$.
 Using this with   \eref{FS2} we derive the lemma for $1<p<\infty$.  

One can derive the lemma in the case $p=1$ and also obtain a constant not depending on $p$ as $p\to 1$ by using a modified Hardy-Littlewood maximal function
$$
M'_0(g)(x):=\sup_{Q\ni x}\left (\frac{1}{|Q|}\int_Q|g(u)|^\mu\, du\right )^{1/\mu}
$$
with $\mu<1$.  The Fefferman-Stein inequality now holds for $p=1$ if this new maximal function is used in place of $M_0$.  For $M_0'$ one
has an analogue of \eref{FS3} where in the first inequality the exponent $d$ is replaced by $d/\mu$.  Thus if $\mu$ is sufficiently close to one so that $\mu(\nu-d)>d$, we again arrive at the lemma for $p=1$.   When $p=\infty$, one can again derive the lemma using the fact that $\nu>2d$ by an analogous argument to the proof of \eref{FS2}.  \endproof
 
 The following is the main result of this section.
\begin{theorem}
\label{nonlinearthm}
Let $1\le p< \infty$ be given, and let $f\in F_{\tau,q}^s$, with
$s \le \k$, $\tau=(1/p+s/d)^{-1}$, and $q:=(1+s/d)^{-1}$.
Then
$$\sig_N(f)_p\le C N^{-s/d}\norm{f}_{F^s_{\tau,q}}.$$
The constant $C$ here is independent of $f$ and $N$.
\end{theorem}

In the proof of the theorem, we will use the following elementary observation.
\begin{lemma}
\label{elemlemma}
Let $\sum_{j=-\infty}^\infty z_j$ be a non-negative series with limit
$Z<\infty$, and  with partial sum sequence $Z_k:=\sum_{j=-\infty}^kz_j$.
For any $\eps>0$, there is a  $C_\eps>0$ depending only on $\eps$ such that 
$$\sum_{j=-\infty}^\infty {z_j\over Z_j^{1-\eps}}\le C_\eps Z^\eps.$$
\end{lemma}

{\bf Proof of the lemma:} For each positive integer $k$, let $j_k$
be the minimal integer for which $Z_{j_k}\ge 2^{-k}Z$, $j_0:=\infty$.
Then
$$\sum_{j=j_{k}}^{j_{k-1}-1}{z_j\over Z_j^{1-\eps}}\le
2^{k(1-\eps)}Z^{\eps-1} 2^{-(k-1)}Z=2Z^\eps 2^{-\eps k}.$$
Summing over all positive $k$, we obtain the stated result.
\endproof

{\bf Proof of the theorem:} Fix $f\in F^s_{\tau,q}$ with $f$ not the zero function, and fix a positive integer 
$N$.   For any given $x\in\Rd$, we consider the set ${\cal V}_x$ of all
$\v\in \calV$ such that $\chi_{I_\v}(x)\neq 0$.   We can order the
$\v\in{\calV}_x$ as follows.  We take any fixed ordering for $E$ and then we
say $\v> \v'$ if either $|\v|>|\v'|$ or $|\v|=|\v'|$ and $e_\v>e_{\v'}$.
Given this order, we now define for each $x\in \Rd$ and   each $\v'\in \calV$, 
$$
M_{q,\v'}(x):= \(\sum_{\v\in\calV_x:\v\ge \v'}
|\v|^{-qs/d}|f_{\v}|^q\chi_{{\iv}}(x)\)^{1/q}.
$$
Notice that  $M_{q,\v'}(x)$  is actually constant on $I_{\v'}$ and so we 
denote  
\beqn
\label{defMq}
M_{q,\v'}:= M_{q,\v'}(x), \quad x\in \ivp.
\eeqn
Also, we clearly have
\beqn
\label{defmqI}
M_{q,\v'}\le M_q(f)(x),\quad x\in \ivp.
\eeqn
We determine our cost function by the rule
$$c_\v=a|\v|^{q}|f_\v|^q M_{q,\v}^{\tau-q},\quad \v\in \calV,$$
where  $a$ will be specified in a moment.
Now $\tau-q\ge 0$ and $q=1-qs/d$.  Therefore, from  (\ref{defmqI}), we see that 
$$c_\v=  \norm{a |\v|^{-qs/d}|f_\v|^q M_{q,\v}^{\tau-q}\chi_{\iv}}_{L_1(\R^d)}\le 
\norm{a |\v|^{-qs/d}|f_\v|^q M_{q}(f)^{\tau-q}\chi_{\iv}}_{L_1(\R^d)}.
$$
Thus,
\begin{eqnarray}
\sum_{\v\in \calV}c_\v&\le&  a \norm{M_q(f)^{\tau-q}
\sum_{\v\in\calV}|\v|^{-qs/d}|f_\v|^q \chi_{\iv}}_{L_1(\Rd)}\cr
&\le&
 a \norm{M_q(f)^{\tau}}_{L_1(\Rd)}= a \norm{f}_{F^s_{\tau,q}}^\tau<\infty,
\end{eqnarray}
where we have used the fact that $\chi_{\iv}\le\chi_{\biv}$ for all $\v\in {\calV}$.
Thus, we can choose $a$ so that $ a\norm{f}_{F^s_{\tau,q}}^\tau=N$, and obtain that $\sum_{\v}c_\v\le N$.

It remains to estimate the $L_p$-error produced by the scheme. In view
of Lemma \ref{maximalfun}, we need to estimate the $L_p$-norm of
$$\sum_{\v\in \calV}[\max(1,c_\v)]^{-\k/d}|f_\v|\chi_{\iv}\le \sum_{\v\in \calV}[\max(1,c_\v)]^{-s/d}|f_\v|\chi_{\iv}\le \sum_{\v\in \calV} c_\v^{-s/d}|f_\v|\chi_{\iv}=:\sum_{\v}E_\v.$$
Here we have used the fact that $s\le k$.

For $x\in \iv$, we have
$$E_\v (x)=a^{-s/d}|\v |^{-qs/d}|f_\v |^{1-qs/d}M_{q,\v }^{(q-\tau)s/d}
=a^{-s/d}|\v |^{-qs/d}|f_\v |^q M_{q,\v }^{\tau/p-q}.$$
Here, we have used the fact that $1-qs/d=q$.
If $p=1$, then $\tau/p-q=0$, and, fixing $x$, we obtain
$$\sum_{  \v\in \calV}E_\v (x)\le  a^{-s/d}M_q^q(f)(x)=a^{-sp/d}M_q^\tau(f)(x).$$
We can prove a similar estimate when $p>1$.  Namely, we fix $x\in\Rd$ and invoke Lemma \ref{elemlemma} with $z_\v :=|\v |^{-qs/d}|f_\v |^q\chi_{\iv}(x)$ using the ordering on ${\calV}_x$ introduced earlier.   Hence,  $M_{q,\v }^q=Z_\v $ and
$M_{q,\v }^{\tau/p-q}=Z_\v ^{\eps-1}$, with $\eps=\tau/(pq)>0$. Also,
$Z\le M_q^q(f)(x)$ again because $\chi_{\iv}\le\chi_{\biv}$. By the lemma,
$$\sum_{x\in \biv }|\v |^{-qs/d}|f_\v |^q M_{q,\v }^{\tau/p-q}\le
C(\tau,p,q)
M_q(f)^{\tau/p}.$$
Thus,
$$(\sum_{\v\in\calV}E_\v(x))^p\le C a^{-sp/d} M_q(f)^\tau(x),$$
and we conclude that, with $A:=\norm{f}_{F^s_{\tau,q}}$
$$\norm{\sum_\v E_\v }_{L_p(\Rd)}\le C a^{-s/d}A^{\tau/p}=
C (a A^\tau)^{-s/d}A=C N^{-s/d}\norm{f}_{F^s_{\tau,q}}.$$
\endproof

\medskip
We can derive from the theorem a corresponding result for the Besov space
$B^s_q(L_\tau)$ (see any of the standard texts for a definition of these
spaces.)  This Besov space is continuously
embedded in $F^s_{\tau,q}$. Hence we obtain

\begin{cor}
Let $1\le p<\infty$ be given, and let $f\in B^s_{q}(L_\tau)$, with
$s \le \k$, $\tau=(1/p+s/d)^{-1}$, and $q:=(1+s/d)^{-1}$.
Then
$$\sig_N(f)_p\le C N^{-s/d}\norm{f}_{B^s_{q}(L_\tau)}.$$
The constant $C$ here is independent of $f$ and $N$.
\end{cor}

Finally,  we compare this theorem with  the classical results on $N$-term 
wavelet approximation given in \cite{DJP}. For wavelet approximation one 
obtains the same bounds with the  assumption 
$f\in B_\tau^s(L_\tau)=F^s_{\tau,\tau}$.  
Since $q<\tau$, the wavelet assumption is (slightly) weaker than what is 
assumed in Theorem \ref{nonlinearthm}.  The two assumptions agree when 
$p=1$.  We do not know if $q$ can be replaced by $\tau$ for other values of 
$p$. We believe that it cannot, i.e., that the value of $q$ in our theorems
is the best possible one.

\medskip
\centerline{\bf Acknowledgment}
We are indebted to Thomas Hangelbroek for  fruitful discussions concerning
this article.

\bibliographystyle{siam}
\bibliography{Amosbib}



\medskip\noindent
Ronald DeVore, Department of Mathematics, University of South Carolina,
Columbia, SC 29208, USA, {\tt devore@math.sc.edu}

\noindent
Amos Ron, Computer Science Department, University of Wisconsin-Madison,
Madison, WI 53706, USA, {\tt amos@cs.wisc.edu}
\end{document}